\documentclass[red,11pt,a4paper]{article}

\usepackage{cite}

\usepackage{amsmath}
\usepackage{amscd}
\usepackage{amssymb}
\usepackage{latexsym}
\usepackage{color}
\usepackage[normalem]{ulem}

 \newcommand{\Deltaa}{\Xi}
 \newcommand{\lt}{\left}
\newcommand{\rt}{\right}
\newcommand{\mms}{\mathcal{S}}

\def\tilde{\widetilde}
\def\hat{\widehat}

\newcommand{\p}{\partial}

\newcommand{\dy}{{\rm d}y}

\newcommand{\pa}{\partial}
\newcommand{\om}{\Omega}

\newtheorem{definition}{Definition}
\newtheorem{theorem}{Theorem}

\newtheorem{corollaire}{Corollary}

\usepackage{./ika}

\title{Pressureless Euler with nonlocal interactions as a singular limit of degenerate Navier-Stokes system}
\author{Jos\'{e} A. Carrillo$^1$, Aneta Wr\'oblewska-Kami\'nska$^2$, Ewelina Zatorska$^3$}
\date{}
\begin{document}
 \maketitle 
\centerline{1. Department of Mathematics, Imperial College London, }
\centerline{London SW7 2AZ, United Kingdom.}
\bigskip
\centerline{2. Institute of Mathematics, Polish Academy of Sciences, }
\centerline{ \'Sniadeckich 8, 00-656 Warszawa, Poland.}
\bigskip
\centerline{3. Department of Mathematics, University College London, }
\centerline{Gower Street, London WC1E 6BT, United Kingdom.}
\bigskip

{\bf Abstract.} We show that weak solutions of degenerate Navier-Stokes equations converge to the strong solutions of the pressureless Euler system with linear drag term, Newtonian repulsion and quadratic confinement. The proof is based on the relative entropy method using the artificial  velocity formulation for the one-dimensional Navier-Stokes system.

\medskip


 \section{Introduction}
Hydrodynamic models of collective behaviour provide a macroscopic description of large groups of interacting individuals. They can be derived from the particle models via BBGKY hierarchies or mean-field limits \cite{hatad, Carrillo2010} and integration of the moments of  the Vlasov-type  kinetic equation seen at the intermediate mesoscopic level. We are particularly interested in the pressureless Euler system with nonlocal interactions modelled by potential $W$ and a linear damping term
 \begin{align}\label{main_eq}
 \begin{cases}
\begin{aligned}
&\pa_t \bar\vr + \pa_x (\bar\vr \bar u) = 0, \cr
&\pa_t (\bar\vr \bar u) + \pa_x (\bar\vr \bar u^2) = -\bar\vr \bar u -(\pa W \ast \bar\vr)\bar\vr, \\
&W(x) = - |x| + \frac{|x|^2}{2},
\end{aligned}
\end{cases}
\end{align}
for $ (t,x) \in \R_+ \times \R$. 
This system is a macroscopic model for individuals whose short-range repulsion is described by the Newtonian  potential $K(x)=-|x|$, the long-range attraction is described by the quadratic confinement $L(x)=\frac12 x^2$, and constant alignment is described by the linear drag term $\bar\vr\bar u$. System \eqref{main_eq}
 was recently considered in \cite{CaChZa} where threshold conditions for the existence of global-in-time solutions emanating from smooth initial data (ref. Theorem 3.1) were derived. It was also proven  that the long-time asymptotic profile   of the density is a  step function determined by the total  mass, the first moment of the initial density, and the initial momentum (ref. Theorem 4.1). 

In this paper we compare the classical solution of \eqref{main_eq} to a weak solution of the corresponding Navier-Stokes type system
\begin{equation}
\begin{cases}
\begin{aligned}
&\pt \vr+\px(\vr u)=0,\\
&\pt\lr{\vr u}+\px(\vr u^2)-\ep\px (\mu(\vr)\px u) + \ep\px p(\vr)
=-\vr u-\vr\px W\ast\vr,\\
&W=-|x|+\frac{|x|^2}{2},
\end{aligned}
\end{cases}
\label{main1}
\end{equation}
for $ (t,x) \in \R_+ \times \R$, where $\ep>0$ denotes a small parameter that will tend to $0$.
The unknowns of system \eqref{main1} are $\vr=\vr(t,x)$ the density and $u=u(t,x)$ the velocity, $p(\vr)=\vr^\gamma$, $\gamma>1$, denotes the barotropic pressure, and $\mu(\vr)=\gamma\vr^\gamma$ denotes the density-dependent viscosity coefficient. Note that the forms of the pressure and the viscosity coefficient are related, we will explain the reason for it later on.  \\
 
\noindent Our main result is the proof that when $\ep\to0$, weak solutions to system \eqref{main1} converge to the strong solution of \eqref{main_eq} as long as the latter exists. The basic idea is to use the momentum equation of \eqref{main1} written in a modified form. Following \cite{HaZa} we introduce an artificial velocity~$v$:
	\eq{\label{defv}
	v = u +  \frac{\ep\gamma}{\gamma-1}\partial_x  \vr^{\gamma-1},
	}
which, at least at the formal level, satisfies the equation
\eq{\label{main_v}
\pt\lr{\vr v}+\px(\vr u v)
=-\vr v-\vr\px W\ast\vr.
}
Let us observe that \eqref{main_v} is very similar to the momentum equation of system \eqref{main_eq}. Moreover, for the limit system \eqref{main_eq} the standard velocity $\bar u$ and the artificial velocity $\bar v$ are in fact the same, since $\ep=0$ in this limit. Using this observation we will construct a relative entropy functional allowing to measure the distance between the weak solutions to the primitive system \eqref{main1} and the strong solutions to the limit system \eqref{main_eq}.

We refer to system \eqref{main1} as the {\emph{degenerate}} Navier-Stokes system because, unlike for usual Navier-Stokes system \cite{Mucha1D},  when the density vanishes the velocity vector field may not be defined. Therefore, the notion of weak solution must be adjusted to deal with this fact. The first attempt to study such system can be found in the work of Veigant and Kazhikhov \cite{VeigantKaz} who considered the initial-value problem on a square for the density dependent viscosity coefficients satisfying additional growth conditions
We will be looking for weak solutions to system \eqref{main1} that satisfy certain energy-entropy estimates. In this context, the first results devoted to weak solutions of the degenerate Navier-Stokes system in the multi-dimensional case are due to  Bresch, Desjardins and coauthors \cite{BD, BDG, BDL, BrDeZa}. In this series of papers they showed essentially the weak sequential stability of solutions to such systems in space dimension $d\geq 1$, and existence of weak solutions to various augmented variants of such system (involving singular pressure, or higher order friction terms). The proof of sequential stability of weak solutions without any of such terms is due to Mellet and Vasseur \cite{MV07}, who combined the entropy estimate, called the Bresch-Desjardins entropy, with additional estimate for the velocity. 
This extra estimate provides sufficient information to prove compactness in the convective term $\vr u\otimes u$. However, construction of weak solutions satisfying all these extra entropy/energy estimates has been a long time unsolved problem. 
Meanwhile, the one-dimensional variant of the problem was addressed by several different authors, and we refer to the work of  Jiang, Xin and Zhang \cite{JXZ} as well as to another work of Mellet and Vasseur \cite{MV07-1D}, where existence of global-in-time strong solutions was proven for various forms of degenerate viscosities. 
Construction of a weak solution in 1D satisfying  all entropy/energy inequalities in the situation when the density may touch the vacuum  was first proven by Li, Li and Xin \cite{LLX} for the initial-boundary problem and then adapted by Jiu and Xin \cite{Jiu} to the whole space case. Very recently, a complete proof of existence in the multi-dimensional case was provided by Vasseur and Yu \cite{VaYu} for shallow water model and by Bresch, Vasseur and Yu \cite{BVY} for the general case. Once again it turned out that the three estimates: classical energy estimate, Bresch-Desjardins entropy estimate and Mellet-Vasseur estimate of the velocity play an essential role in the proof. 
As far as we know, the existence of solutions to \eqref{main1} with the nonlocal interaction forces is so far restricted to the Navier-Stokes-Poisson type systems \cite{Ding,LiuYuan}. 
In this paper we show that the three major estimates necessary to repeat the construction from \cite{VaYu} or from \cite{Jiu} are true for the Navier-Stokes system with more general interaction terms, provided they are combined with the estimates of higher moments of the density. With these a-priori estimates at hand, the compactness arguments from \cite{MV07}, see also \cite{HaZa, Zat1, Zat2}, allow to show the sequential stability of weak solutions
to \eqref{main1}. This is shown in the Section \ref{Sec:3} of the paper. Later on, in Section \ref{Sec:4}, we show that when $\ep\to0$ the weak solutions to system \eqref{main1} coincide with the strong solutions of system \eqref{main_eq} as long as the latter exist. This is achieved thanks to the relative entropy inequality that relies on the reformulation of the momentum equation  \eqref{main_v}. Similar form of the entropy inequality was used by Haspot in \cite{Haspot1D} to prove the weak-strong uniqueness of the solutions to the degenerate  Navier-Stokes equations in one dimension and in \cite{BrNoVi2017} for applications to various  singular limits problems.
It is also worth to mention a result of Brenier \cite{Brenier} who used the monotonic rearrangement operator to reduce the degenerate one-dimensional Navier-Stokes system (also the Navier-Stokes-Poisson system) on a torus  to an elementary differential equation with noise. In this work the pressure and viscosity are linked to each other through a given smooth and strictly convex function, and the limit passage $\ep\to 0$ does not rely on the relative entropy functional.

\section{Preliminaries and the main result}
In this section we introduce the concept of weak solution to system \eqref{main1}, the strong solution to system \eqref{main_eq}, and formulate our main result.
\subsection{Weak solutions to the Navier-Stokes system}
We start from introducing the hypothesis on the initial data
\eq{\label{initiald}
(\vr(t,\cdot), u(t,\cdot))|_{t=0} = ( \vr_0, u_0),}
for which we assume that
\begin{equation}
\begin{aligned}
&\vr_0\geq 0, \quad \vr_0\in L^1(\R)\cap L^\infty(\R), \quad |x|^{\kappa+2}\vr_0\in L^1(\R),\quad
(\vr_0^{\gamma-\frac{1}{2}})_x\in L^2(\R),  \\
&\lim_{|x|\to\infty}\vr_0(x) u_0(x)=0,\quad \vr_0 u_0^2\in L^1(\R),\quad \vr_0|u_0|^{2+\kappa}\in L^1(\R),
\end{aligned}
\label{ini}
\end{equation}
where $0<\kappa\leq\min\{2\gamma-1,\frac2\gamma\}$. The total initial mass and momentum are given s.t.
\begin{equation}\label{int_01_1}
\begin{split}
& M_0 := \intO{\vr_0(x)} \qquad\mbox{and}
\qquad
M_1=\intO{\vr_0(x) u_0(x)},
\\
& 0<M_0<\infty,\quad |M_1|<\infty.
\end{split}
\end{equation}

\begin{definition}\label{Def:1}
    For  fixed $\ep>0$ the pair of functions $(\vr,\sqrt{\vr}u)$ is called a weak solution to the system \eqref{main1} with initial data \eqref{initiald} satisfying \eqref{ini} and \eqref{int_01_1} if:
        \begin{equation}\label{weak_vr_u}
        \begin{split}
        &\limsup_{|x|\to+\infty}|\vr(t,x) u(t,x)|=0,\quad for\ a.a. \ t\in(0,T),\\
        &  \vr\in C_{w} ([0,T]; L^\gamma(\R)), \ \vr\in L^\infty(0,T; L^\gamma(\R)\cap L^1(\R)), \\
         &
         \sqrt{\vr} u \in L^\infty(0,T;L^2(\R)),\\ 
       & \partial_x (\vr^{\gamma-\frac{1}{2}})\in L^\infty(0,T; L^2(\R)),\\
       & \vr^{\gamma-\frac{1}{2}} \in L^\infty(0,T;L^\infty(\R)).
        \end{split}       
        \end{equation}
Moreover, the following weak formulation of the continuity equation is satisfied 
\begin{equation}
\int_{\R}\vr\psi(t_2)\,\dx-\int_{\R}\vr\psi(t_1)\,\dx=
\int^{t_2}_{t_1}\!\!\!\int_{\R}(\vr\p_t\psi+\vr u\p_x\psi)\,\dxdt
\label{2.3}
\end{equation}
 for any $T\geq t_2\geq t_1\geq 0$ and any $\psi\in C^1( [t_1,t_2]\times \R)$, and denoting $\vr v=\vr u+\ep\px\vr^\gamma$ the following equality holds
\begin{equation}
\int_{\R}\vr\psi(t_2)\dx-\int_{\R}\vr\psi(t_1)\,\dx=
\int^{t_2}_{t_1}\!\!\!\int_{\R}\lr{\vr\p_t\psi+\vr v\p_x\psi-\ep\px\vr^\gamma\px\psi}\dxdt.
\label{new2.3}
\end{equation}
Furthermore, the the weak formulation of the momentum equation
\eq{
\int_{\R}\vr_0u_0\psi(0)\dx+\intTO{\lr{(\vr u)\p_t\psi+(\vr u^2+\ep\vr^\gamma)\p_x\psi}}{-\ep\langle\mu(\vr)\p_x u,\p_x\psi\rangle}\\
=\intTO{\vr u\psi}+\intTO{\vr(\px W\ast\vr) \psi},
\label{2.4}
}
is satisfied for any $\psi\in C^{\infty}_{c}([0,T)\times\R)$, where the diffusion term is defined as follows:
\eq{
&\langle\mu(\vr)\p_x u,\p_x\psi\rangle\\
&=-\gamma\intTO{\vr^{\gamma-\frac{1}{2}}\sqrt{\vr}u{\p_{xx}}\psi}-
\frac{2\gamma}{2\gamma-1}\intTO{\p_x\lr{\vr^{\gamma-\frac{1}{2}}}\sqrt{\vr} u{\px\psi}}.
\label{2.5}
}

\end{definition}
\begin{rmk}
Let us emphasize that all terms in the above formulation are well defined. In particular, notice that 
    $$ 
    \intTO{\lr{(\vr u)\p_t\psi+(\vr u^2}} =
    \intTO{\lr{\sqrt{\vr}(\sqrt{\vr}u)\p_t\psi+((\sqrt{\vr}u)^2}}.$$
\end{rmk}
With this basic definition at hand, following Dafermos \cite{dafermos1979} and Haspot \cite{Haspot1D}, we introduce the relative entropy (energy) functional 
\begin{equation}\label{relE1}
\begin{split}
&{\cal E}(\vr,v|\tilde\vr,\tilde v)(t):=\\
&\quad 
\int_\R \left(
\frac{\vr (v-\tilde v)^2}{2}+\frac{1}{2} (\vr-\tilde\vr)W\ast(\vr-\tilde\vr)
+\ep\lr{H(\vr) - H'(\tilde\vr)(\vr - \tilde\vr) - H(\tilde\vr)}\right)\dx,
\end{split}
\end{equation}
where $v$ is given by \eqref{defv}, $ H(s) =  \frac{1}{\gamma} s^\gamma$, and $\tilde\vr(t,x),\tilde v(t,x)$ are smooth functions defined on $[0,T]\times\R$, such that
\eq{\label{general_f}\tilde\vr>0, \ on\ [0,T]\times\R,\quad \tilde\vr(t,x)\tilde v(t,x)\to 0 \ as\ x\to\pm\infty.}
To simplify presentation, we immediately assume  that $\tilde\vr(t,x),\tilde v(t,x)$ satisfy the limit system \eqref{main_eq}, so in fact we take
$$\tilde\vr(t,x)=\bar\vr(t,x),\quad \tilde v(t,x)=\bar u(t,x).$$
{Note that since $H(s)$ is a strictly convex function we have
$$\ep\lr{H(\vr) - H'(\bar\vr)(\vr - \bar\vr) - H(\bar\vr)}\geq \ep H''(\xi)R^2=\ep(\gamma-1)\xi^{\gamma-2}R^2,$$
for some $\xi$ strictly between $\vr$ and $\bar\vr$. 
 Therefore, this term gives nonnegaive contribution in  \eqref{relE1}.
Concerning the non-local terms in the relative entropy functional \eqref{relE1}, the same form was obtained in \cite{CaFaGw} where the weak-strong uniqueness of solutions to the compressible Euler system with nonlocal terms in multi-dimensional setting was studied, and in \cite{BrezinaMacha} where a viscous approximation of measure-valued solutions was considered. Because the potential $W$ is symmetric, the nonlocal terms give positive contrigution in \eqref{relE1}. Moreover, the Newtonian repulsion term $K(x)=-|x|$ implies in particular that the $W^{-1,2}$ distance between the solutions $\vr$ and $\bar\vr$ is controlled. }

\noindent Following some tedious calculations presented in Section \ref{Sec:4}, we will derive the relative entropy inequality at the level of sufficiently smooth approximation of system \eqref{main1}:
\eq{\label{rel_ineq}
&[{\cal E}(\vr,v|\bar\vr,\bar u)]_{t=0}^{t=\tau}+\lr{\frac{\ep\gamma}{\gamma-1}}^2\int_0^\tau\!\!\intR{\vr\lr{\partial_x \varrho^{\gamma-1}}^2}\,\dt+2\ep\int_0^\tau\!\!\intO{\vr^\gamma\bar\vr}\,\dt\\
&\hspace{2cm}\leq \int_0^\tau {\cal R}(\vr,v,\bar\vr,\bar u)\,\dt,}
for almost all $\tau\in(0,T)$, with the reminder term given by
\eq{\label{reminder}
& {\cal R}(\vr,v,\bar\vr,\bar u)=-\intR{\vr (v-\bar u)(u-\bar u)(\px\bar v+1)}
-\intR{\vr (u-\bar u)\ep\frac{\gamma}{\gamma-1}\partial_x \bar \varrho^{\gamma-1}}\\
&\qquad\qquad\qquad+\intR{(\vr-\bar\vr)\bar v\px W\ast(\vr-\bar\vr)}
+2\ep\intO{\vr^\gamma\bar\vr}-\gamma\intO{\px \bar v F},
}
and 
$$F= F(\vr,\bar\vr)=\ep\lr{H(\vr) - H'(\bar\vr)(\vr - \bar\vr) - H(\bar\vr)}.$$
\begin{rmk}
A weak solution to \eqref{main1} that satisfies inequality similar to \eqref{rel_ineq} but for general functions as in \eqref{general_f} is sometimes referred to as a suitable weak solution.  It is not clear, however, whether this inequality could be justified for the weak solutions. Therefore, we derive it at the level of approximation, estimate the reminder, and pass to the limit with the parameter of the approximation in the resulting formula.
\end{rmk}

\subsection{Strong solutions to the Euler system}\label{Sec:strong}
As mentioned in the introduction, the target system is pressureless compressible Euler-type of system with nonlocal interactions \eqref{main_eq}. The critical thresholds for the initial data for existence of global in time strong solutions were determined in \cite{CaChZa}.
 Following their approach we consider system \eqref{main_eq} supplemented by the initial values for the density and the velocity
 \eq{\label{ini_main_eq}
(\bar\vr(t,\cdot), \bar u(t,\cdot))|_{t=0} = (\bar \vr_0,\bar u_0)\in C^2(\Omega_0)\times C^3(\Omega_0),
 }
 s.t. $\bar\vr_0>0$ on $\Omega_0$, 
where  $\om_0:=\om(0)$ is either an open interval $(a_0,b_0)$ (in which case $\bar\vr(t, x)$ is extended by $0$ outside $\Omega(t)$) or the whole line $\om_0=\R$. 
{In addition, we assume boundedness of total initial mass and momentum (we assume that both are the same as in the previous system, see \eqref{int_01_1}, for relative entropy inequality)}
\begin{equation}\label{int_01_2}
M_0 = \int_{\Omega_0}\bar \vr_0(x)\, \dx \qquad\mbox{and}
\qquad
M_1=\int_{\om_0} \bar \vr_0(x)\bar u_0(x)\,\dx,
\end{equation}
as well as boundedness of the first moment of the density and integrability of the initial momentum
\eq{\label{Rem_31}
\int_\R |x| \bar \vr_0(x)\,dx <\infty \qquad \mbox{and} \qquad \int_\R \bar \vr_0(x) |\bar  u_0(x)|\,dx <\infty.
}
\begin{definition}\label{Def:2}
We say that a pair of functions $(\bar \vr,\bar u)$ is a classical local-in-time solution to \eqref{main_eq} with the initial data \eqref{ini_main_eq}, if:
\begin{itemize}
\item There exists  time $T>0$ such that $\bar\vr(t,x)>0$, and  $\bar \vr,\bar u$ are $C^1$ and $C^2$ respectively in the set $\{(t,x)\in [0,T) \times \Omega(t)\}$
and  $\bar \vr$ and $\bar u$ satisfy  the equations \eqref{main_eq} pointwisely in $\{(t,x)\in [0,T) \times \Omega(t)\}$ with initial data \eqref{ini_main_eq}.
     \item The characteristics $\eta(t,x)$ associated to $\bar u$ defined by 
\eq{\label{eq_tra}
\frac{d \eta(t,x)}{dt} =  \bar u(t,\eta(t,x)) \quad \mbox{with} \quad \eta(0,x) = x \in \om_0
}
 are diffeomorphisms for all $t\in [0,T)$ with $\Omega(t)=\eta(t, \Omega_0)$.
\end{itemize}
\end{definition}
For the sake of completeness let us recall  Theorem 3.1 from \cite{CaChZa} stating the global existence or finite-time blow up of classical solutions to \eqref{main_eq} depending  on the initial data $\px \bar u_0$, $\bar \vr_0$ and  $M_0$.
\begin{theorem}[Carrillo, Choi, Zatorska,  2016]
\label{Th:mainCCZ} 
Assume that $(\vr,u)$ is a classical solution to the system \eqref{main_eq} with initial data \eqref{ini_main_eq}, then:

\noindent {\bf Case A:} If $1 - 4M_0 > 0$, the solution blows up in finite time if and only if there exists a $x^* \in \om_0$ such that
\[
\pa_x \bar u_0(x^*) < 0, \quad M_0 - 2\bar \vr_0(x^*) < \lambda_1 \pa_x \bar u_0(x^*),
\]
and
\[
2\bar \vr_0(x^*) \leq (\lambda_1 \pa_x \bar u_0(x^*) - M_0 + 2\bar \vr_0(x^*))^{-\lambda_2/\sqrt{\Deltaa}}(\lambda_2 \pa_x \bar u_0(x^*) - M_0 + 2\bar \vr_0(x^*))^{\lambda_1/\sqrt{\Deltaa}},
\]
where
$$\lambda_1 := \frac{-1 + \sqrt{1 - 4M_0}}{2}, \quad \lambda_2 := \frac{-1 - \sqrt{1 - 4M_0}}{2}, $$
$$
\Deltaa := 1 - 4M_0\,.
$$
{\bf Case B:} If $1 - 4M_0 = 0$, the solution blows up in finite time if and only if there exists a $x^* \in \om_0$ such that
\[
\pa_x \bar u_0(x^*) < \min\lt\{0, 4\bar \vr_0(x^*) - \frac12\rt\},
\]
and
\[
\ln \lt(\frac{8\bar \vr_0(x^*)}{8\bar \vr_0(x^*) - 2\pa_x \bar u_0(x^*) - 1}\rt) \leq \frac{2\pa_x \bar u_0(x^*)}{8\bar \vr_0(x^*) - 2\pa_x \bar u_0(x^*) - 1}.
\]
{\bf Case C:} If $1 - 4M_0 < 0$, the solution blows up in finite time if and only if there exists a $x^* \in \mms_1 \cup \mms_2\cup\mms_3\cup\mms_4$ where 
\begin{align}\label{condi_c3}
\begin{aligned}
\mms_1 \!\!&:=\!\! \lt\{ x \in \om_0\!:\!\pa_x C_5(x) < 0, \pa_x C_6(x) >0, \frac{2\bar \vr_0(x)}{M_0} -  C_7(x)\exp\!\lt( \frac{C_8(x)}{\sqrt\square}\rt) \!\!\leq 0 \rt\},\cr
\mms_2 \!\!&:=\!\! \lt\{ x \in \om_0\!:\! \pa_x C_5(x) > 0, \pa_x C_6(x) <0, \frac{2\bar \vr_0(x)}{M_0} - C_7(x) \exp\!\lt( \frac{C_8(x) - \pi}{\sqrt\square}\rt) \!\!\leq 0 \rt\},\cr
\mms_3 \!\!&:= \!\!\lt\{ x \in \om_0\!:\! \pa_x C_5(x) < 0, \pa_x C_6(x) <0, \frac{2\bar \vr_0(x)}{M_0} - C_7(x) \exp\!\lt( \frac{C_8(x) - \pi}{\sqrt\square}\rt) \!\!\leq 0 \rt\},\cr
\mms_4 \!\!&:=\!\! \lt\{ x \in \om_0\!:\! \pa_x C_5(x) > 0, \pa_x C_6(x) >0, \frac{2\bar \vr_0(x)}{M_0} - C_7(x) \exp\!\lt( \frac{C_8(x) - 2\pi}{\sqrt\square}\rt) \!\!\leq 0 \rt\}
\end{aligned}
\end{align}
where
$$ \square := -\Deltaa \,,$$
\[
\pa_x C_5(x) = \pa_x \bar u_0(x), \quad \pa_x C_6(x) = \frac{2}{\sqrt\square} \lt(-\frac12 \pa_x \bar u_0(x) - M_0(x) + 2\bar \vr_0(x)\rt),
\]
\[
C_7(x) := \lt( \frac{2\sqrt\square}{1 + \square}\rt)\sqrt{\frac{1 + \square}{\square}(\pa_x \bar u_0(x))^2 + \frac{4}{\square}\lt(2\bar \vr_0(x) - M_0\rt)\lt(2\bar \vr_0(x) - M_0 - \pa_x \bar u_0(x)\rt)},
\]
and
\[
C_8(x) := \arctan \lt( \frac{\sqrt{\square} \pa_x\bar  u_0(x)}{4\bar \vr_0(x) - 2M_0 - \pa_x \bar u_0(x)}\rt).
\]
Moreover, for all cases, if there is no finite-time blow-up, then the classical solution $(\bar \vr,\bar u)$ exists globally in time.
\end{theorem}

The proof of this theorem was based on derivation of an explicit formula for the solution $u=u(t,\eta(t,x))$ in the Lagrangian coordinates. 
This allowed, for instance, to determine the exact form of $\eta(t,x)$ from \eqref{eq_tra}, and from there the exact form of the density on the characteristics, using the continuity equation
\eq{\label{bar0}
\bar\vr(t,\eta(t,x))\partial_x \eta(t,x)=\bar \vr_0(x),\quad x\in\Omega_0.}

For the purposes of this paper it is important to observe the following estimates of the classical solutions to \eqref{main_eq} for large $x$.

\begin{corollaire}\label{co_1}
Let $\Omega_0 = \R$ and  let
\begin{subequations}\label{ini_add}
\eq{\label{ini_add1}
|\bar u_0(x)|\leq C(1+|x|),\ \px \bar u_0, \ \px^2 \bar u_0 \in L^\infty (\R),}
\eq{\label{ini_add2}
{\lim_{|x|\to\infty}\eta_x(t,x)\geq c>0},
}
\end{subequations}
for any $t\in(0,\infty)$.\\
Then the solution  $\bar{u}$ to \eqref{main_eq} has the following properties:
\begin{equation}\label{EP_u_x}
| \bar{u}(t,x) | \leq c_1+ c_2|x| \quad \mbox{ for } t\in (0,T), \ x\in \R,
\end{equation}
\begin{equation}\label{EP_pxu_c}
| \partial_x \bar{u}(t,x) | \leq c_3 \quad \mbox{ for } t\in (0,T), \ x\in \R,
\end{equation}
were $c_i$, $i=1,2,3$ are positive constants depending on initial data and $T$. 

\end{corollaire}

\pf Let us denote $\bar{w}(t,x) := \bar{u} (t,\eta(t,x))$. Following \cite[Section~2]{CaChZa} we have
\begin{equation}\label{eq_diff}
\pa_{tt} \bar{w} + \pa_{t} \bar{w} + M_0 \bar{w} = M_1 e^{-t}, \quad t > 0, \quad \bar{w}_0 = \bar{u}_0
\end{equation}
and the initial data $\pa_{t}\bar{w} (t,x)\big|_{t=0} =\bar{w}'_{0}(x)$ are given through 
\begin{equation}\label{vt0}
\bar{w}'_{0}(x)= -\bar{w}_0(x) - (x+1) M_0 + \int_{\R} y \bar\vr_0(y)\dy + 2\int_{-\infty}^x \bar\vr_0(y) \dy  \quad \mbox{for} \quad x \in \R.
\end{equation}
Moreover, we can compute 
\begin{equation}\label{epe}
\eta(t,x) = x + \int_0^t \bar{w}(s,x){\rm \,d}s \quad \mbox{and} \quad \pa_x \eta(t,x) = 1 + \int_0^t \pa_x \bar{w}(s,x){\rm \,d}s.
\end{equation}
Depending on the size of the initial mass $M_0$, as long as the solution exists, it satisfies:\\
$\bullet$ {\bf Case A} ($1 > 4M_0$):
\begin{equation}\label{sol_formA}
\bar{w}(t,x) =  C_1 e^{\lambda_1  t} + C_2 e^{\lambda_2 t} + \frac{M_1}{M_0} e^{-t} ,
\end{equation}
$\bullet$ {\bf Case B} ($1 =4M_0$):
\begin{equation}\label{sol_formB}
\bar{w}(t,x) =C_3 e^{-t/2} + C_4 t \,e^{-t/2} + \frac{M_1}{M_0} e^{-t},
\end{equation}
$\bullet$ {\bf Case C} ($1 <4M_0$):
\begin{equation}\label{sol_formC}
\bar{w}(t,x) =C_5 e^{-t/2} \cos\lt( \frac{\sqrt{4M_0 - 1}}{2}t\rt) + C_6 e^{-t/2} \sin\lt( \frac{\sqrt{4M_0 - 1}}{2}t\rt) + \frac{M_1}{M_0} e^{-t},
\end{equation}
where $\lambda_1$, $\lambda_2$, and $C_i,i=1,\cdots,6$ are given by
\begin{subequations}\label{coeff}
\begin{align}
\lambda_1 &:= \frac{-1 + \sqrt{1 - 4M_0}}{2}, \quad \lambda_2 := \frac{-1 - \sqrt{1 - 4M_0}}{2}, \label{coeff_l}\\
C_1 &:= \frac{1}{\lambda_2 - \lambda_1}\lt( \lambda_2 \bar{w}_0 - \bar{w}_0' + \lambda_1\frac{M_1}{M_0}\rt), \quad 
C_2 := \frac{1}{\lambda_2 - \lambda_1}\lt( -\lambda_1 \bar{w}_0 + \bar{w}_0' - \lambda_2\frac{M_1}{M_0}\rt), \label{coeff_12}\\
 C_3 &:= \bar{w}_0 - \frac{M_1}{M_0}, \quad C_4:= \frac{\bar{w}_0}{2} + \bar{w}_0' + \frac{M_1}{2M_0},\label{coeff_34}\\
C_5&:= \bar{w}_0 - \frac{M_1}{M_0}, \quad \mbox{and} \quad C_6 = \frac{2}{\sqrt{4M_0 - 1}}\lt(\bar{w}_0' + \frac{\bar{w}_0}{2} + \frac{M_1}{2M_0}\rt).\label{coeff_56}
\end{align}
\end{subequations}
For abbreviation, we set 
$$
\Deltaa := 1 - 4M_0\qquad \mbox{and} \qquad \square := -\Deltaa \,.
$$
{Note that due to this explicit formula for the solution, assumption \eqref{ini_add2} is in fact an assumption for the initial condition $\bar u_0, \bar u_0'$. Indeed, since in all the cases A, B, C recalled above, the sollution $\bar w$ is expressed as combination of functions that are integrable in time ($e^{-ct},\ te^{-ct}, \ e^{-c_1t}\cos(c_2t)$ for some positive constants $c,c_1,c_2$), the assumption \eqref{ini_add2} is met for example in the case when the $L^\infty$ norms of $\px \bar w_0$ and $\px\bar w_0'$ are sufficiently small.}

From these representations we also immediately check that  that $\bar{w}$ and $\partial_x \bar{w}$ satisfy:\\
 $\bullet$ {\bf Case A}
 \begin{equation}\label{cAwpxw}
 \begin{split}
&  |\bar{w}(t,x) | \leq C_{A,1} e^{\max(\lambda_1,\lambda_2)t} (|x| + 1) \quad\mbox{ as }|x| \to \infty,\\
& |\partial_x \bar{w} (t,x) | \leq C_{A,2} e^{\max(\lambda_1,\lambda_2)t} \quad\mbox{ as }|x| \to \infty,
 \end{split}
 \end{equation}
 $\bullet$ {\bf Case B}
 \begin{equation}\label{cBwpxw}
 \begin{split}
 &  |\bar{w}(t,x) | \leq C_{B,1} (1+t)) e^{- \frac{t}{2}}  (|x| + 1)  \quad\mbox{ as }|x| \to \infty,\\
& |\partial_x \bar{w} (t,x) | \leq C_{B,2} e^{- \frac{t}{2}}  \quad\mbox{ as }|x| \to \infty,
 \end{split}
 \end{equation}
 $\bullet$ {\bf Case C}
 \begin{equation}\label{cCwpxw}
 \begin{split}
 &  |\bar{w}(t,x) | \leq C_{C,1} e^{- \frac{t}{2}}  (|x| + 1)\quad\mbox{ as }|x| \to \infty,\\
& |\partial_x \bar{w} (t,x) | \leq C_{C,2} e^{- \frac{t}{2}} \quad\mbox{ as }|x| \to \infty. 
 \end{split}
 \end{equation}
One can summarise that for any fixed $T$ we have that
$$|\bar{w}(t,x)| \leq c|x| \quad \mbox{ as } |x| \to \infty.$$
In order to come back to Eulerian coordinates we notice that 
	$$
	| \bar u (t,\eta(t,x))| \leq c |x|  \quad \mbox{ as } |x| \to \infty
	$$
and consequently 
\begin{equation}\label{bueta-1}
| \bar u (t, y )| \leq c | \eta^{-1} (t,y)|  \quad\mbox{ with } y= \eta(t,x) \quad \mbox{ as } |x| \to \infty
\end{equation}
Note that $\px\eta|_{t=0}=1$, moreover from the second equality in \eqref{epe} and formulas for $\px\bar w$ in all three cases we verify that for classical solutions $\lim_{t\to\infty}\px\eta(t,x)=1$. The positive bound from below for $\px\eta$ for $|x|\to \infty$ follows also from \eqref{epe} together with the assumption  \eqref{ini_add2}. Since for classical solutions $\px\eta$ is a continuous function, its minimal value is bounded away from 0. The bound from above follows directly from the assumption \eqref{ini_add1}, hence we have
\eq{\label{good_eta}
0< \underline{c} \leq \partial_x \eta \leq \overline{c} < \infty.} 
Consequently, the growth of $\eta^{-1}$ is not faster than linear. Next writing 
$$\partial_y \bar{u}(t,y) = \partial_x \bar{w}(t, \eta^{-1}(t,y)) \partial_y(\eta^{-1})(t,y), $$
since by \eqref{cAwpxw}, \eqref{cBwpxw}, \eqref{cCwpxw}, we in fact have that $| \partial_y \bar{w} | \leq C $ for any fixed time, 
we immediately find that  
$$|\bar{u} (t,y) - \bar{u} (t,0)| \leq c|y|,$$
where $\bar{u}(t,0) = \bar{w}(t,\eta^{-1}(t,0))$. \hfill $\Box$

Additionally to what has been said above let us observe that $\bar\vr$ behaves similarly to $\bar{\vr}_0$. 
Indeed, one can write 
\eq{\label{rhorho0}
\bar{\vr}(t,\eta(t,x)) \partial_x \eta(t,x) = \bar{\vr}_0}
and consequently 
$$\underline{c} \leq \left| \frac{\bar\vr(t, \eta(t,x))}{\bar{\vr}_0(x)} \right| \leq \overline{c}. $$
In particular, if $\bar{\vr}_0$ is positive on $\R$, then $\bar{\vr}$ remains positive on $\R$ for all $t\in (0,T)$. 

\subsection{The main result}
The main theorem of this paper consists of two parts.
We first prove that the weak solution to \eqref{main1} is sequentially stable  for $\ep>0$ fixed. This means that having a smooth approximation of system \eqref{main1} satisfying all energy-entropy bounds, there exists a subsequence whose limit satisfies the Definition \ref{Def:1} together with the relative entropy inequality. Secondly, we prove the convergence of the weak solutions to \eqref{main1}, as $\ep\to0$, to the strong solution of \eqref{main_eq}, emanating from the same initial data, as long as the latter exist.
\begin{theorem}\label{Thm:main} Assume that $\gamma\in(1,3/2]$.\\
(i) Let $(\vr_n,u_n)$ be a sequence of weak solutions  to \eqref{main1} in the sense of Definition \ref{Def:1} satisfying energy inequalities  \eqref{energye}, \eqref{BD_ineq}, and \eqref{testMV3}, where $v$ is defined by \eqref{defv}. Let the initial conditions
\eq{
\vr_n|_{t=0}=\vr_0,\quad u_n|_{t=0}= u_0
}
satisfy assumptions \eqref{ini} and \eqref{int_01_1}.  Then, up to a subsequence $(\vr_n,\sqrt{\vr_n}u_n)$ converges strongly to $(\vr_\ep,\sqrt{\vre}u_\ep)$  a weak solution of \eqref{main1} satisfying energy inequalities 
{\eqref{energye}, \eqref{BD_ineq}, and \eqref{testMV3}}. More precisely
\eq{\label{c_delta}
&\vr_n\to\vre \quad\text{strongly\ in } C([0,T]\times\R_{loc}),\\
&\px\lr{\vr_n^{\gamma-{\frac12}}}\to\px\lr{\vre^{\gamma-{\frac12}}}\quad\text{weakly in } L^2(0,T; L^2(\Omega)),\\
&\sqrt{\vr_n} u_n\to\sqrt{\vr} u_\ep, \quad\text{strongly in } L^{2}(0,T; L^{2}(\R)),\\
&\vr_n v_n\to\vr v_\ep, \quad\text{strongly in } L^{2}(0,T; L^{2}(\R)),\\
&\gamma\vr_n^\gamma\px u_n\to\Lambda_\ep \quad\text{weakly in } L^2(0,T; L^2_{loc}(\R)),
}
where $\Lambda_\ep$ satisfies
\begin{equation}
\intTO{\Lambda_\ep \phi}=-\gamma\intTO{\vre^{\gamma-\frac{1}{2}}\sqrt{\vre}u_\ep\p_x\phi}-\frac{2\gamma}{2\gamma-1}\intTO{\p_x(\vre^{\gamma-\frac{1}{2}})\sqrt{\vre}u_\ep \phi}
\label{2.9}
\end{equation}
for any $\phi\in C^\infty_c((0,T)\times\R)$.

\medskip

(ii) 
{Let $(\vr_n,u_n)$ satisfy in addition the entropy inequality \eqref{rel_ineq}-\eqref{reminder}}, with $(\bar\vr,\bar u)$ being the strong solution to \eqref{main_eq} on the time interval $(0,T)$  in the sense specified in Definition \ref{Def:2}. 
 Let the initial data $(\bar\vr_0,\bar u_0)=(\vr_0,u_0)$ satisfy \eqref{ini} and \eqref{int_01_1} together with
 \begin{equation}
(\vr_0^{\gamma-1})_x\in L^2(\R).
\label{ini_lim}
\end{equation}
Then the limiting weak solution obtained in \eqref{c_delta} converges to the strong solution of \eqref{main_eq} in the following sense
$${\rm ess}\sup_{t\in(0,T)}{\cal E}(\vre,v_\ep|\bar\vr,\bar u)(t)\to 0\quad as\ \ep\to0.$$
\end{theorem}

\begin{rmk}
Note that the initial conditions for the primitive as well as the limit problem are the same. This assumption could be relaxed, but we do not focus on this aspect here. Note, in particular, that assumptions  \eqref{ini} and \eqref{int_01_1} imply \eqref{Rem_31} that is necessary for existence of strong solutions to the limit system. 
\end{rmk}
\begin{rmk}
Thanks to  \eqref{ini} and \eqref{int_01_1}, via \eqref{good_eta} and \eqref{rhorho0} we can verify that 
 $\bar\vr\bar u\to 0$ and $x^{3/2}\bar\vr\to 0$ for $x\to\infty$ and that $\bar\vr\in L^{\gamma+1}((0,T)\times\R)$. Thanks to \eqref{ini_lim} we also get that $\bar\vr\in L^{\gamma+1}((0,T)\times\R)$. These properties of the limit solution are needed in order to close the relative entropy estimate. 
\end{rmk}

The rest of this paper is devoted to the proof of this theorem. The proof of existence of approximation $(\vr_n, u_n)$  could be obtained following \cite{Jiu} or \cite{VaYu} and is left for the future study.

\section{A-priori estimates for the Navier-Stokes system}\label{Sec:3}
The purpose of this section is to provide various essential a-priori estimates for the approximate solutions $(\vr_n, u_n)$ to the primitive system \eqref{main1}.
We will show that these a-priori estimates  provide  uniform with respect to $n$ bounds that allow to deduce the weak sequential stability of solutions. Similar reasoning was performed before in \cite{MV07, HaZa}.

In what follows we drop the subindex $n$ when no confusion can arise,  and we assume that the assumptions of Theorem \ref{Thm:main} are satisfied. In particular, $(\vr,u)$ are sufficiently smooth functions satisfying equations of system \eqref{main1} pointwise, s.t. $\vr\geq0$, and $\lim_{|x|\to\infty}\vr(t,x)u(t,x)\to0$. We clearly do not expect that $u(t,x)$ decays at infinity itself, but we a-priori assume that $\px u$ is bounded at infinity, in accordance with what is known for the limiting system, see \eqref{EP_pxu_c}.
\subsection{Conservation of mass and momentum}
It is straightforward to deduce that $\vr$ is nonnegative.  Integrating the continuity equation with respect to space variable, and using decay of $\vr(t,x) u(t,x)$ at infinity we obtain that
\eq{
\frac{ \rm d}{\dt}\intO{ \vr}=0,}
and so the total mass is conserved, in particular
\eq{\label{mass_cons}
\|\vr_n\|_{L^\infty(0,T; L^1(\R))}\leq C,
}
uniformly w.r.t. $n$.

Similarly, integrating the momentum equation, we check that since $\px W=-sgn(x)+2x$ is antisymmetric, the nonlocal term integrates to zero, and hence we have
\eq{
\frac{ \rm d}{\dt}\intO{ \vr u}=-\intO{\vr u}.}
Integrating with respect to time, we find
\eq{\label{mom}
\intO{ \vr u (t)}= e^{-t}\intO{ \vr_0 u_0}=e^{-t} M_1,}
and so the total momentum is bounded if the initial momentum  is.
\subsection{The basic energy estimate}
We recall the classical energy inequality that is derived by multiplying the momentum equation of \eqref{main1} by $u$ and integration by parts. For the nonlocal term, this leads to
\eq{\label{nonle}
-\int_\R\vr u\px W\ast\vr \dx=\int_\R\px(\vr u) W\ast\vr  \dx =-\int_\R\pt\vr W\ast\vr \dx=- \frac{\rm d}{\dt}\int_\R \frac{1}{2} \vr W\ast\vr \dx, 
}
where we used the continuity equation to get the second equality.
Similarly, for the pressure term we have that
$$\ep\int_\R \partial_x \varrho^\gamma u \dx 
=  \ep\frac{1}{\gamma-1} \frac{ \rm d}{\dt} \intR{\vr^{\gamma}}.$$
Consequently, the energy inequality reads
\eq{\label{energye}
\frac{\rm d}{\dt} \intO{ \left(\frac12 \varrho u^2 + \frac{\ep}{\gamma-1}\vr^{\gamma} 
+ \frac{1}{2} \varrho W \ast \varrho
\right)}
+ \ep\gamma \intO{ \varrho^\gamma (\partial_x u)^2} + \intO{ \varrho u^2 } \leq0.
}
Note that $W=-|x|+x^2\geq -C$, so the l.h.s. can be made positive, by adding to $W$ a constant, and by using the fact that the mass is conserved.

Integrating \eqref{energye} over time interval $[0,T]$ we obtain the following estimates 
\eq{
&\|\sqrt{\vr_n}u_n\|_{L^\infty(0,T; L^2(\R))}\leq C,\\
&\|\vr_n\|_{L^\infty(0,T; L^\gamma(\R))}\leq C,\\
&\|\vr_n^{\gamma/2}\px u_n\|_{L^2(0,T; L^2(\R))}\leq C
}
 are uniform w.r.t. $n$. Smoothness of $\vr_n$ and $u_n$ for $n$ fixed together with the bound from above implies also that $\lim_{|x|\to0}\sqrt{\vr_n(t,x)}u_n(t,x)=0$, which will become important later on.
\subsection{Bresch-Desjardins entropy estimate}
Testing \eqref{main_v} by solution $v$, and using the continuity equation we obtain
	$$
	\frac{\rm d}{\dt}\intO{ \frac{\varrho v^2}{2} } + \intO{ \varrho v^2 } 
	+ \intO{ \varrho \lr{u + \ep \frac{\gamma}{\gamma-1}\partial_x \varrho^{\gamma-1} }\partial_x W \ast \varrho }= 0,
	$$
where in the last term we  used the definition of $v$	\eqref{defv}. Proceeding as in \eqref{nonle} we therefore get
	$$
	\frac{\rm d}{\dt}\intO{ \left(  \frac{1}{2} \varrho v^2 
	+ \frac{1}{2} \varrho W \ast \varrho  \right) }
	+ \intO{ \varrho v^2 }
	+ \ep\intO{ \partial_x \varrho^{\gamma}  \partial_x W \ast \varrho } = 0.
	$$
Let us now focus on the last term on l.h.s., integrating by parts we get
\eq{\label{noi}
\ep \intO{ \partial_x \varrho^{\gamma}  \partial_x W \ast \varrho } =-\ep\intO{\vr^\gamma \partial^2_x W\ast\vr},
}
where we used the fact that for $n$ fixed  $\vr$ is smooth and integrable, therefore $\vr_n(t,x)\to 0$ when $|x|\to+\infty$.
Due to definition of $W$ we have $\partial^2_x W = -2 \delta + 1$, and so
$$
-\ep\intO{\vr^\gamma \partial^2_x W\ast\vr}= 
2\ep  \intO{ \varrho^{\gamma+1}}- \ep \intR{\vr}\intR{\vr^\gamma}=2\ep  \intO{ \varrho^{\gamma+1}}- \ep M_0\intR{\vr^\gamma},
$$
which gives rise to the mathematical entropy inequality
\eq{\label{BD_ineq}
\frac{\rm d}{\dt}\intO{ \left(  \frac{1}{2} \varrho v^2 
	+ \frac{1}{2} \varrho W \ast \varrho \dx \right) }
	+ \intO{ \varrho v^2 }
	+ 2\ep  \intO{ \varrho^{\gamma+1}}\leq  \ep M_0\intR{\vr^\gamma}.
}
After integration over time we have 
	\eq{\label{BD_main}
	\intOB{\frac{\varrho v^2}{2} + \frac{1}{2}   \varrho W \ast \varrho (t)} +
	\intTOB{ \varrho v^2 + 2\ep  \varrho^{\gamma+1}} \leq C + \ep CT.
	}
 From this estimate it follows in particular by \eqref{defv} that uniformly w.r.t. $n$ we have
\begin{equation}\label{rhonormg}
\ep\|\partial_x(\vr_n^{\gamma-\frac{1}{2}})\|_{L^\infty(0,T; L^2(\R))}\leq C,
\end{equation}
for $\ep$ small enough and so, by the Sobolev imbedding, since $\vr \in L^1\cap L^\gamma(\R)$,
\begin{equation}\label{rhonorm}
\ep\|\vr_n^{\gamma-\frac{1}{2}}\|_{L^\infty(0,T; L^\infty(\R))}\leq C.
\end{equation}

\subsection{Higher moments estimates}
Let us first check the behaviour of the center of mass of the density. Multiplying continuity equation by $x$ and integrating by parts we obtain
\begin{equation}\label{fm1}
\Dt\int_\R {x\vr} \dx = \int_\R {\vr u} \dx.
\end{equation}
 Indeed, this can be justified using as a test function for the continuity  equation $x\Psi_k$ where
$\Psi_k \in C^{\infty}_c(\R)$, $\Psi_k = 1$ for $x \in [-k+1, k-1]$, $\Psi_k = 0$ for $x \in (-\infty,- 2k -1) \cup (2k+1, \infty)$, $| \partial_x \Psi_k | \leq \frac{1}{k}$.
Then 
$$\frac{\rm d}{\dt} \int_\R x \vr \Psi_k \dx = - \int_\R \partial_x (\vr u ) x \Psi_k \dx
= \int_{\R} \vr u \Psi_k \dx + \int_{\R} \vr u x \partial_x \Psi_k \dx.
$$
Since for a.a. $t\in(0,\infty)$, $\vr u = \sqrt{\vr} \sqrt{\vr}u \in L^1(\R)$ (notice $\sqrt{\vr} \in L^2(\R)$ and $\sqrt{\vr}u \in L^2(\R)$), we may pass with $k \to \infty$ and obtain \eqref{fm1}.
 
Next, using the momentum estimate \eqref{mom}, we obtain that
\begin{equation}\label{M1t}
M_2(t):=\int_\R {x\vr(t)}\dx = \int_\R{x\vr_0} \dx +(1-e^{-t})M_1.
\end{equation}
With this at hand, we go one level higher and we estimate the second moment of the density. Using again the continuity equation we have
\eq{\label{x2rho}
\Dt\intO{x^2\vr}=2\intO{x \vr u}.}
This again can be justified using appropriate test function in the continuity equation. We now take
$\Phi_k \in C^{\infty}_c(\R)$, such that $\Phi_k = 1$ for $x \in [-k^2+1, k^2-1]$, $\Phi_k = 0$ for $x \in (-\infty,- k^2 -k-1) \cup (k^2+k+1, \infty)$, $| \partial_x \Phi_k | \leq \frac{1}{k^2}$.
Then integrating by parts
$$\frac{\rm d}{\dt} \int_\R x^2 \vr \Phi_k \dx = - \int_\R \partial_x (\vr u ) x^2 \Phi_k \dx
= \int_{\R} x \vr u \Phi_k \dx + \int_{\R} \vr u x^2 \partial_x \Phi_k \dx.
$$
Since we already know that  $x \vr u \in L^1(\R)$, we may pass with $k \to \infty$ and obtain \eqref{x2rho}.

The r.h.s. of \eqref{x2rho} can be evaluated using the momentum equation, we have
\eq{\label{mmom}
&\Dt\intO{x \vr u}\\
&\quad=\intO{\vr u^2}-\ep\gamma\intO{\vr^\gamma\px u}+\ep\intO{\vr^\gamma}-\intO{x\vr u}-\intO{x\vr\px W\ast\vr},
}
hence we get
\begin{equation*}
	\intO{ x (\varrho  u)(t) } = e^{-t}\intO{ x \vr_0 u_0 }
	+ \int_0^t e^{s-t} f(s) {\rm \,d}s,
	\end{equation*}
where
\eqh{
	f(s) =  \intO{\lr{\varrho u^2  -  \ep\gamma\varrho^\gamma \partial_x u +  \ep \varrho^\gamma-x\vr\px W\ast\vr
	}(s)}.
}
From the energy estimate we can check that $f\in L^2((0,T))$. Clearly,   $\| \varrho u^2 \|_{L^\infty(0,T; L^1(\Omega))} \leq C$, $\sqrt{\ep} \| {\varrho}^{\frac{\gamma}{2}} \|_{L^{\infty}(0,T; L^2(\Omega))} \leq C$, $\sqrt{\ep} \| {\varrho}^{\frac{\gamma}{2}} \partial_x u \|_{L^2(0,T;L^2(\Omega))} \leq C$. The only issue is the last term, for which we can write
\eqh{
&-\intO{x\vr\px W\ast\vr}=-\int_\R\intO{x\vr(x)\partial W(x-y)\vr(y)}\,\dy\\
&=-\int_\R\intO{(x-y)\vr(x)\partial W(x-y)\vr(y)}\,\dy-\int_\R\intO{y\vr(x)\partial W(x-y)\vr(y)}\,\dy\\
&=-\int_\R\intO{(x-y)\vr(x)\partial W(x-y)\vr(y)}\,\dy-\int_\R\intO{x\vr(y)\partial W(y-x)\vr(x)}\,\dy\\
&=-\int_\R\intO{(x-y)\vr(x)\partial W(x-y)\vr(y)}\,\dy+\int_\R\intO{x\vr(y)\partial W(x-y)\vr(x)}\,\dy.
}
Thus, it follows that
\begin{equation}\label{term25}
-\intO{x\vr\px W\ast\vr}=-\frac12\intO{\vr (x\px W)\ast\vr}.
\end{equation}
Now note that  $|x\px W|\leq C(W+1)$, then we have thanks to the energy estimate  \eqref{energye} and conservation of mass that
the term \eqref{term25} is bounded uniformly in time.
Consequently we deduce that 
$$ \intO{ x (\varrho  u)(t) } \leq C\lr{T, \intO{ x \vr_0 u_0 }},$$
and so, by \eqref{x2rho} we infer that 
	\begin{equation}\label{second_mom}
	\intO{ x^2 \varrho (t) }  \leq c \left( T, \intO{ x^2 \varrho_0 }, \intO{ x \vr_0 u_0 } \right).
	\end{equation}
Finally, we want  to estimate higher order moments. Due to continuity equation for any $\kappa>0$ we have
\eq{\label{xkrho}
\Dt\intO{x^2|x|^\kappa\vr}=(2+\kappa)\intO{x|x|^\kappa \vr u}.}
We cannot proceed as previously using the momentum equation. However, applying the Young inequality to the r.h.s of \eqref{xkrho} we can write
\eq{\label{xkrho2}
\Dt\intO{|x|^{\kappa+2}\vr}&\leq C\intO{|x|^{\kappa+1} \vr |u|}\\
&\leq C\intO{\lr{|x|^{\kappa+2}\vr}^{\frac{\kappa+1}{\kappa+2}}\vr^{\frac{1}{\kappa+2}}|u|}\\
&\leq C_1(\kappa) \intO{|x|^{\kappa+2}\vr}+C_2(\kappa)\intO{\vr|u|^{2+\kappa}}.
}
And that will be used later in forthcoming section (see \eqref{testMV3}).

\subsection{Mellet-Vasseur velocity estimate}
The final estimate is the improved estimate of the velocity a'la Mellet and Vasseur \cite{MV07}.  It is obtained by testing the momentum equation of the system \eqref{main1} by $u |u|^\kappa$ for $0<\kappa\leq \min\{2\gamma-1,\frac{2}{\gamma}\}$, and by testing the continuity equation by $\frac{|u|^{2+\kappa}}{2+\kappa}$. Summing up the obtained expressions we obtain
\eq{\label{testMV}
&\Dt\intO{\frac{\vr|u|^{2+\kappa}}{2+\kappa}}+\ep\gamma(\kappa+1)\intO{\vr^\gamma|u|^\kappa|\px u|^2}+\intO{ \vr |u|^{2+ \kappa}} \\
&\hspace{2cm}=-\ep\intO{\px\vr^\gamma u |u|^\kappa}-\intO{\vr u|u|^\kappa\px W\ast\vr}=I_1+I_2 .
} 
Due to decay of $\vr$ and $\vr u^2$  at infinity the lack of boundary term coming from integration by parts in the viscosity term is justified provided $\kappa\leq2\gamma-1$. Under the same restriction
 $I_1$ can be integrated by parts and estimated as follow
\eqh{I_1&=-\ep\intO{\px\vr^\gamma u|u|^\kappa}
=\ep(\kappa+1)\intO{\vr^\gamma|u|^\kappa\px u }\\
&\leq
\frac{\ep(\kappa+1)}{2}\intO{\vr^\gamma|u|^\kappa|\px u|^2}+ \frac{\ep(\kappa+1)}{2}\intO{\vr^\gamma|u|^\kappa}.
}
The first term can be absorbed by the l.h.s. of \eqref{testMV}, while to control the second term, we use Young's inequality to write for $\kappa\leq2$
\eqh{ \frac{\ep(\kappa+1)}{2}\intO{\vr^\gamma|u|^\kappa}= \frac{\ep(\kappa+1)}{2}\intO{\vr^{\gamma-\frac{\kappa}{2}}\vr^{\frac{\kappa}{2}}|u|^\kappa}\\
\leq C(\kappa)\intR{\vr u^2}+ \ep C(\kappa)\intR{\vr^{\frac{2\gamma-\kappa}{2-\kappa}}}
}
for  $\ep$ small enough and so, the first term can be controlled by the energy estimate \eqref{energye}, while the second one is bounded by the $L^1$ norm of $\vr^{\gamma+1}$ from \eqref{BD_main} provided $\kappa\leq\frac{2}{\gamma}$.
For the estimate of $I_2$, we use the explicit form of $\px W\ast \vr$, i.e.
\eq{
\px W\ast \vr=M_0-2\int_{-\infty}^x \vr(t,y)\,\dy+xM_0-\intO{x\vr}.
}
This means that 
\eq{
I_2&\leq\left|
-\intO{\vr u|u|^\kappa\px W\ast\vr}\right|\\
&\leq 3M_0\intO{\vr |u|^{\kappa+1}}+M_0\intO{|x|\vr |u|^{\kappa+1}}
+\intO{|x|\vr}\intO{\vr |u|^{\kappa+1}}.
}
Using the estimates of the second moment of $\vr$ \eqref{second_mom} and of the total mass we find that
$$\intO{ |x| \vr  (t)} \leq c(M_0, m_0, T, \intO{ x^2 \vr_0 }),$$
and therefore 
\eqh{
I_2 &\leq C\intO{|x|\vr |u|^{\kappa+1}}
+C\intO{\vr |u|^{\kappa+1}}\\
&=C\intO{\lr{|x|^{2+\kappa}\vr}^{\frac{1}{2+\kappa}} \lr{\vr^{\frac{1}{2+\kappa}}|u|}^{\kappa+1}}+
C\intO{\vr^{\frac{1}{2+\kappa}} \lr{\vr^{\frac{1}{2+\kappa}}|u|}^{\kappa+1}}.
} 
Now, we use the  Young inequality again with $p=\kappa+2$, $p'=\frac{\kappa+2}{\kappa+1}$, similarly as in \eqref{xkrho2}, for both terms at the same time, so that we get
\eq{
I_2
\leq C(\kappa)\lr{ \intO{|x|^{2+\kappa}\vr}+\intO{\vr}+\intO{\vr|u|^{\kappa+2}}}.
}
Summarizing, the formula \eqref{testMV} might be now rewritten as:
\eq{\label{testMV2}
\Dt\intO{\frac{\vr|u|^{2+\kappa}}{2+\kappa}}+ \frac{\ep(\kappa+1)(2\gamma -1)}{2}\intO{\vr^\gamma |u|^\kappa|\px u|^2}+\intO{ \vr |u|^{2+ \kappa} } \\
\leq C(\kappa)\lr{1+ \intO{|x|^{2+\kappa}\vr}+\intO{\vr|u|^{\kappa+2}}}.}
In the above $C(\kappa)$ depends also on $M_0,$ $m_0,$ $T,$ $\int_\Omega x^2 \vr_0 \dx$.
In order to deduce some useful estimates, we now add it to formula \eqref{xkrho2}, to get
\eq{\label{testMV3}
\Dt\lr{\intO{\frac{\vr|u|^{2+\kappa}}{2+\kappa}}+\intO{|x|^{2+\kappa}\vr}}+ \frac{\ep(\kappa+1)(2\gamma -1)}{2}\intO{\vr^\gamma |u|^\kappa|\px u|^2}\\
\leq C(\kappa)\lr{1+ \intO{|x|^{2+\kappa}\vr}+\intO{\vr|u|^{\kappa+2}}}}
for $\kappa$ small enough, 
and we conclude by Gronwall's inequality that
if $T<\infty$, $\vr_0 |u_0|^{2+\kappa} \in L^1(\R)$, $|x|^{2+\kappa}  \vr_0 \in L^1(\R)$ then   
	\begin{equation}\label{higher_m}
	\begin{split}
	&\| \vr_n |u_n|^{2+\kappa}\|_{L^\infty(0,T; L^1(\R))} \leq C,\\
	&  \||x|^{2+\kappa}  \vr_n \|_{L^\infty(0,T; L^1(\R))}  \leq C,
	\end{split}
	\end{equation}
where $C$ depends on: $\kappa$,  $M_0,$ $m_0,$ $T,$ $\int_\Omega x^2 \vr_0 \dx$, $\vr_0 |u_0|^{2+\kappa}$, $|x|^{2+\kappa}  \vr_0$, but is uniform w.r.t. $n$.\\
We summarize the findings of this section in the following lemma.
\begin{lemma}
Let  $(\vr_n, u_n)$ be a sequence of sufficiently smooth approximate solutions to system \eqref{main1} with the initial data \eqref{initiald} satisfying \eqref{ini} and \eqref{int_01_1}. Then, uniformly in $n$ we have
\eq{
&\|\vr_n\|_{L^\infty(0,T; L^1(\R))}\leq C,\\
&\|\sqrt{\vr_n}u_n\|_{L^\infty(0,T; L^2(\R))}\leq C,\\
&\|\vr_n\|_{L^\infty(0,T; L^\gamma(\R))}\leq C,\\
&\|\vr_n^{\gamma/2}\px u_n\|_{L^2(0,T; L^2(\R))}\leq C,\\
&\ep\|\partial_x(\vr_n^{\gamma-\frac{1}{2}})\|_{L^\infty(0,T; L^2(\R))}\leq C,\\
&\ep\|\vr_n^{\gamma-\frac{1}{2}}\|_{L^\infty(0,T; L^\infty(\R))}\leq C\\
&\| \vr_n |u_n|^{2+\kappa}\|_{L^\infty(0,T; L^1(\R))} \leq C,\\
&  \||x|^{2+\kappa}  \vr_n \|_{L^\infty(0,T; L^1(\R))}  \leq C
}
for $\kappa$ small enough.
\end{lemma}

\section{Relative entropy} \label{Sec:4}
In this section we first derive the relative entropy inequality and then the relative entropy estimate. We again assume that $(\vr_n, u_n)$ are sufficiently smooth functions satisfying equations of system \eqref{main1} pointwise, s.t. $\vr\geq0$, and $\lim_{|x|\to\infty}\vr(t,x)u(t,x)\to0$.

Recall that due to \eqref{defv}, at the level of sufficiently regular solutions, the  system  \eqref{main1} is equivalent to the following one:
\begin{equation}
\begin{cases}
\begin{aligned}
&\pt \vr+\px(\vr u)=0,\\
&\pt\lr{\vr v}+\px(\vr uv)
=-\vr v-\vr\px W\ast\vr.
\end{aligned}
\end{cases}
\label{main_aa}
\end{equation}
Setting $\bar v=\bar u$ in the Euler system \eqref{main_eq} we obtain the following equations:
\begin{equation}
\begin{cases}
\begin{aligned}
&\pt \bar\vr+\px(\bar\vr \bar u)=0,\\
&\pt\lr{\bar\vr \bar v}+\px(\bar\vr {\bar u}{\bar v})=-\bar\vr \bar v-\bar\vr\px W\ast\bar\vr.
\end{aligned}
\end{cases}
\label{main_bb}
\end{equation}
Note that in both cases $W=-|x|+\frac{|x|^2}{2}$.
\subsection{Derivation of the relative entropy inequality}
We will use the equations \eqref{main_aa} and \eqref{main_bb} to construct the relative entropy functional  similar to the one obtained in \cite{Haspot1D}, modulo the nonlocal attraction/repulsion term and linear damping term. Note also that in \cite{Haspot1D} it was assumed that $\bar v$ is a bounded function which is not the case here. Therefore, all integrations by parts and estimates must take this fact into account.
In this part of the proof it is  important that both systems \eqref{main_aa} and \eqref{main_bb} are considered on the whole domain $\R$. Adapting our arguments to the free boundary problem would require some tedious extensions of the solution to the limit problem preserving the Sobolev norms, and we leave this for future research.

For brevity of formulas, we introduce the notation:
$$V=v-\bar v,\quad U=u-\bar u,\quad R=\vr-\bar\vr,$$
and so, subtracting equations of system \eqref{main_bb} from equations of system  \eqref{main_aa}, respectively, we obtain:
\begin{equation}
\begin{cases}
\begin{aligned}
&\pt R+\px(\vr U)+\px(R\bar u)=0,\\
&(\vr\pt+\vr u\px)V+\vr V=-R(\pt\bar v+\bar u\px\bar v)-R\bar v-\vr U\px\bar v-\vr\px W\ast\vr+\bar \vr\px W\ast\bar\vr.
\end{aligned}
\end{cases}
\label{merge}
\end{equation}
Next, we  reduce the right hand side of the new momentum equation using the fact that in the momentum equation of \eqref{main_bb} all terms are multiplied by $\bar\vr$. 
Note that the classical solution to the limit system has positive density in whole domain provided this was the case for the initial condition. Therefore, by  the continuity equation of \eqref{main_bb}, the momentum equation may be reduced to
$$
\pt\bar v+\bar u\px\bar v=-\bar v- \partial_x W\ast\bar\vr .
$$
So, inserting it to the momentum equation of \eqref{merge}, we get
\begin{equation}
\begin{cases}
\begin{aligned}
&\pt R+\px(\vr U)+\px(R\bar u)=0,\\
&(\vr\pt+\vr u\px)V+\vr V=-\vr U\px\bar v -\vr\px W\ast\vr+ \vr\px W\ast\bar\vr \ .
\end{aligned}
\end{cases}
\label{merge_a}
\end{equation}
Our task now is to derive an analogue of the classical energy estimate for the system \eqref{merge_a}. We start from multiplying  the second equation of \eqref{merge_a} by $V$ and integrating over the whole space to get

\begin{equation}\label{rel_1}
\Dt\int_R {\frac{\vr V^2}{2}} = - \int_\R{\vr V^2} \dx  -\int_\R {\vr V U\px\bar v} \dx -\int_\R{\vr\lr{\px W\ast(\vr-\bar\vr)}V}\dx,
\end{equation}
where we also used the continuity equation of \eqref{main_aa} tested by $\frac{1}{2} V^2$ to obtain the energy term on the l.h.s of \eqref{rel_1}.
Although the first  term on  the r.h.s. of \eqref{rel_1} has already a right sign we further transform it using $V  = U + \ep\frac{\gamma}{\gamma-1}\partial_x \varrho^{\gamma-1}$:
\eq{\intR{\vr V^2}=& \intR{\vr VU}+\intR{\vr V\ep\frac{\gamma}{\gamma-1}\partial_x \varrho^{\gamma-1}}\\
=&\intR{\vr VU}+\intR{\vr U\ep\frac{\gamma}{\gamma-1}\partial_x \varrho^{\gamma-1}}+\intR{\vr\lr{\frac{\ep\gamma}{\gamma-1}}^2\lr{\partial_x \varrho^{\gamma-1}}^2}\\
=&\intR{\vr VU}+\intR{\vr U\ep\frac{\gamma}{\gamma-1}\partial_x \bar \varrho^{\gamma-1}} \\
&+\intR{\vr U\ep\frac{\gamma}{\gamma-1}\px\lr{ \varrho^{\gamma-1}- \bar\varrho^{\gamma-1}}}+\intR{\vr\lr{\frac{\ep\gamma}{\gamma-1}}^2\lr{\partial_x \varrho^{\gamma-1}}^2}.
}
Inserting this back to \eqref{rel_1} and rearranging the terms we get
\eq{\label{rel_2a}
&\Dt\intR{\frac{\vr V^2}{2}}+\intR{\vr U\ep\frac{\gamma}{\gamma-1}\px\lr{ \varrho^{\gamma-1}- \bar\varrho^{\gamma-1}}}+\intR{\vr\lr{\frac{\ep\gamma}{\gamma-1}}^2\lr{\partial_x \varrho^{\gamma-1}}^2}\\
&=-\intR{\vr V U(\px\bar v+1)}
-\intR{\vr U\ep\frac{\gamma}{\gamma-1}\partial_x \bar \varrho^{\gamma-1}}
-\intR{\vr\lr{\px W\ast R}V}.
}
In the next step we will complement the kinetic part of the relative entropy appearing on the l.h.s. of \eqref{rel_2a} by the potential part that is linked to the non-local forces. To this purpose we evaluate
\eq{\label{h-1}
&\Dt\intR{\frac{1}{2} (\vr-\bar\vr)W\ast(\vr-\bar\vr)}=\intR{(\vr-\bar\vr)W\ast(\pt\vr-\pt\bar\vr)}\\
&=-\intR{(\vr-\bar\vr)W\ast\px(\vr u-\bar\vr\bar u)}=\intR{(\vr u-\bar\vr\bar u)\px W\ast(\vr-\bar\vr)}.
}
Note that above we only used the fact that $W$ is symmetric, and not the integration by parts. Substituting $u=v-\ep\frac{\gamma}{\gamma-1}\px\vr^{\gamma-1}$, $R=\vr-\bar\vr$, and $\bar u=\bar v$ we obtain
\eq{\label{h-2}
&\Dt\intR{\frac{1}{2} RW\ast R}
=\intR{(\vr v-\ep\vr \frac{\gamma}{\gamma -1}\px\vr^{\gamma-1}-\bar\vr\bar v)\px W\ast R}\\
&=\intR{R \bar v\px W\ast R}+\intR{\vr(v-\bar v)\px W\ast R}
-
\ep\intR{\px\vr^\gamma\px W\ast R}\\
&=\intR{R \bar v\px W\ast R}+\intR{\vr V\px W\ast R}
-
\ep\intR{\px\vr^\gamma\px W\ast R}.
}
Summing up  \eqref{rel_2a} and \eqref{h-2} we obtain
\eq{\label{rel_2b}
&\Dt\intR{\frac{\vr V^2}{2}}+
\Dt\intR{\frac{1}{2} RW\ast R}\\
&\quad+\intR{\vr U\ep\frac{\gamma}{\gamma-1}\px\lr{ \varrho^{\gamma-1}- \bar\varrho^{\gamma-1}}}+\lr{\frac{\ep\gamma}{\gamma-1}}^2\intR{\vr\lr{\partial_x \varrho^{\gamma-1}}^2}\\
&=-\intR{\vr V U(\px\bar v+1)}
-\intR{\vr U\ep\frac{\gamma}{\gamma-1}\partial_x \bar \varrho^{\gamma-1}}\\
&\quad+\intR{R\bar v\px W\ast R}
-
\ep\intR{\px\vr^\gamma\px W\ast R}.
}
To deal with the last term on the r.h.s.  we will separate $W(x)$ into two parts $K(x)=-|x|$, and $L(x)=\frac{x^2}{2}$.  Note that $ \px L\ast(\vr-\bar\vr)  = 0 $, indeed we have
\begin{equation}\label{Lzero}
\begin{split}
\px L\ast(\vr-\bar\vr)  &=\int_\R(x-y)(\vr(t,y)-\bar\vr(t,y)) \dy\\
&=x\int_\R(\vr(t,y)-\bar\vr(t,y))\ \dy-\int_\R (x\vr-x\bar\vr) \dx\\
& = -\int_\R(x\vr-x\bar\vr)\dx=0.
\end{split}
\end{equation}
The two last equalities follow from the fact the total masses for both $\vr$ and $\bar\vr$ as well as their first moments are the same. To see the equality between the moments it is enough to integrate the momentum equation for the target system \eqref{main_eq}, and to check that it implies the same formula  \eqref{M1t} as for the $\ep$-dependent system, and so $\intR{x\bar\vr}=\intR{x\vr}$.

Using \eqref{Lzero} we immediately simplify the last term in \eqref{rel_2b} as follows
\eq{\label{good_bad}
&-\ep\int_\R{\px\vr^\gamma \px W\ast R}\dx= -\ep\int_\R{\px\vr^\gamma \px K\ast R} \dx
=\ep\intR{\vr^\gamma \partial_{xx}K\ast R}\\
&=-2\ep\intR{\vr^\gamma(\vr-\bar\vr)}=-2\ep\intO{\vr^{\gamma+1}}+2\ep\intO{\vr^\gamma\bar\vr},
}
where we used the explicit form of $K$ to compute
$\partial_{xx} K(x)=-2\delta_0(x),$ and the fact that $\vr^\gamma(t,x)\to 0$ as $|x|\to+\infty$ and $\px  K\ast R$ is bounded, so the boundary term in the integration by parts disappears.
Note that the first term on the r.h.s. of \eqref{good_bad} has a good sign  and so it can be moved to the l.h.s. of \eqref{rel_2b}.

Let us now focus on the  the third term on the l.h.s. of  \eqref{rel_2b} and show that it has the right meaning of the distance. Indeed, integrating by parts and using the continuity equation from \eqref{merge}, we get
\eq{\label{rel_3}
\intR{\vr U\ep\frac{\gamma}{\gamma-1}\px\lr{ \varrho^{\gamma-1}- \bar\varrho^{\gamma-1}}}&=-\intR{\px(\vr U)\ep\frac{\gamma}{\gamma-1}\lr{ \varrho^{\gamma-1}- \bar\varrho^{\gamma-1}}}\\
&=\intR{(\pt R+\px(R\bar u))\ep\frac{\gamma}{\gamma-1}\lr{ \varrho^{\gamma-1}- \bar\varrho^{\gamma-1}}}.
}
Note that the boundary terms $[\vr u(\vr^{\gamma-1}-\bar\vr^{\gamma-1})]_{x\to-\infty}^{x\to+\infty}$ dissapears because $\vr u\to 0$ for $|x|\to +\infty$, and $[\vr \bar u(\vr^{\gamma-1}-\bar\vr^{\gamma-1})]_{x\to-\infty}^{x\to+\infty}$ disappears because for $|x|\to+\infty$, we know that $\bar u$ grows at most linearly \eqref{EP_u_x}, while $x\vr$ tends to 0 due to \eqref{second_mom}, and assumption that $\vr$ is smooth, and $(\vr^{\gamma-1}-\bar\vr^{\gamma-1})$ remains bounded. 
Now we define our free energy $F(R,\bar\vr)$ as follows
\eq{
F(R,\bar\vr)=\frac{\ep}{\gamma}(R+\bar\vr)^\gamma-\ep\bar\vr^{\gamma-1}R-\frac{\ep}{\gamma}\bar\vr^\gamma.
}
Note that because $R+\bar\vr=\vr$, we have that
\eq{\label{FR}
\frac{\partial F}{\partial R}=\ep((R+\bar\vr)^{\gamma-1}-\bar\vr^{\gamma-1})=\ep(\vr^{\gamma-1}-\bar\vr^{\gamma-1}),
}
later on we will also need
\eq{\label{Frho}
\frac{\partial F}{\partial \bar\vr}=\ep\lr{(R+\bar\vr)^{\gamma-1}-(\gamma-1)R\bar\vr^{\gamma-2}-\bar\vr^{\gamma-1}}.
}
Using \eqref{FR} in \eqref{rel_3} we get that
\eq{\label{rel_4}
&\intR{\vr U\ep\frac{\gamma}{\gamma-1}\px\lr{ \varrho^{\gamma-1}- \bar\varrho^{\gamma-1}}}=\intR{(\pt R+\px(R\bar u))\frac{\gamma}{\gamma-1}\frac{\partial F}{\partial R}}\\
&=\frac{\gamma}{\gamma-1}\intR{\pt R \frac{\partial F}{\partial R}}+\frac{\gamma}{\gamma-1}\intR{\bar u\px R \frac{\partial F}{\partial R}}+\frac{\gamma}{\gamma-1}\intR{\px\bar u R \frac{\partial F}{\partial R}}\\
&=\frac{\gamma}{\gamma-1}\lr{\intR{\pt F}-\intR{\pt \bar\vr \frac{\partial F}{\partial \bar\vr}}+\intR{\bar u\px F}-\intR{\bar u\px \bar\vr \frac{\partial F}{\partial \bar\vr}} +\intR{\px\bar u R \frac{\partial F}{\partial R}}},
}
where to pass to the last line we used the fact that $\partial_s F(R,\bar\vr)=\partial_s R \frac{\partial F}{\partial R}+\partial_s \bar\vr \frac{\partial F}{\partial \bar\vr}$, for $s=t,x$.

Next, we eliminate the second and fourth terms in above using the continuity equation for $\bar\vr$, from which it follows that:
\eq{-\pt\bar\vr-\bar u\px\bar\vr=\bar\vr \px\bar u ,}
therefore, integrating the third term on the r.h.s. of \eqref{rel_4} by parts we find that 
\eq{\label{rel_5}
&\intR{\vr U\ep\frac{\gamma}{\gamma-1}\px\lr{ \varrho^{\gamma-1}- \bar\varrho^{\gamma-1}}}\\
&=\Dt\intR{\frac{\gamma}{\gamma-1}F}+\frac{\gamma}{\gamma-1}\intR{\px\bar u \lr{-F+\bar\vr\frac{\partial F}{\partial \bar\vr}+ R \frac{\partial F}{\partial R}}}.
}
Indeed, let us notice that the integration by parts is justified here, similarrly as for estimates for higher order moments, by Corollary~\ref{co_1} and the fact that we keep center of mass and second moment for both densities bounded due to assumptions.
Finally, using  \eqref{FR} and \eqref{Frho}, we find that
$$  \bar\vr\frac{\partial F}{\partial \bar\vr}+ R \frac{\partial F}{\partial R}=\gamma F .$$
After this simplification \eqref{rel_2b} can be rewritten in its final form as
\eq{\label{rel_6}
&\Dt\int_\R 
\left( 
\frac{\vr V^2}{2}+\frac{1}{2} RW\ast R+ \frac{\gamma}{\gamma-1}F
\right) \dx\\
&\quad+\lr{\frac{\ep\gamma}{\gamma-1}}^2\intR{\vr\lr{\partial_x \varrho^{\gamma-1}}^2}+2\ep\intO{\vr^{\gamma+1}}\\
&\leq-\intR{\vr V U(\px\bar v+1)}
-\intR{\vr U\ep\frac{\gamma}{\gamma-1}\partial_x \bar \varrho^{\gamma-1}}\\
&\quad+\intR{R\bar v\px W\ast R}
+2\ep\intO{\vr^\gamma\bar\vr}-\gamma\intO{\px \bar u F}.
}
Note that we have
\begin{equation}\label{relE}
\begin{split}
&\int_\R 
\left( 
\frac{\vr V^2}{2}+\frac{1}{2} RW\ast R+ \frac{\gamma}{\gamma-1}F
\right) \dx\\
&=\int_\R \left(
\frac{\vr (v-\bar v)^2}{2}+\frac{1}{2} (\vr-\bar\vr)W\ast(\vr-\bar\vr)
+\ep\lr{H(\vr) - H'(\bar\vr)(\vr - \bar\vr) - H(\bar\vr)}\right)\dx\\
&={\cal E}(\vr,v|\bar\vr,\bar u)(t)
\end{split}
\end{equation}
for $ H(s) =  \frac{1}{\gamma} s^\gamma$, as specified in \eqref{rel_ineq}, and therefore, after integration of \eqref{rel_6} over $(0,\tau)$, we obtain \eqref{rel_ineq} with the reminder as in \eqref{reminder}.

\subsection{Estimate of the reminder}
In what follows we estimate the terms on the r.h.s. of the relative entropy inequality \eqref{rel_6}, i.e.
\eq{\label{rel_66}
{\cal R}(\vr,v,\bar\vr,\bar u)&=-\intR{\vr V U(\px\bar v+1)}
-\intR{\vr U\ep\frac{\gamma}{\gamma-1}\partial_x \bar \varrho^{\gamma-1}}\\
&\quad+\intR{R\bar v\px W\ast R}
+2\ep\intO{\vr^\gamma\bar\vr}-\gamma\intO{\px \bar u F}\\
&=I_1+I_2+I_3+I_4+I_5.
}

 To this purpose we  use the general a-priori estimates derived in Section \ref{Sec:3} as well as the l.h.s. of \eqref{rel_6}.

\medskip
\noindent{\emph{Estimate of $I_1$}.} For the first term, we use the fact that
$$U=u-\bar u=v-\ep\frac{\gamma}{\gamma-1}\partial_x \varrho^{\gamma-1}-\bar v=V-\ep\frac{\gamma}{\gamma-1}\partial_x \varrho^{\gamma-1},$$
therefore
\begin{equation}\label{I1}
\begin{split}
I_1 & =-\int_\R \vr V U(\px\bar v+1)\dx\\
&=
-\int_\R \vr V^2(\px\bar v+1) \dx +
\int_\R \vr V\ep\frac{\gamma}{\gamma-1}\partial_x \varrho^{\gamma-1}(\px\bar v+1)\dx \\
&\leq 
-\int_\R
\vr V^2(\px\bar v+1) \dx 
+\frac{1}{4}
\int_\R
\lr{\frac{\ep\gamma}{\gamma-1}}^2\vr(\px\vr^{\gamma-1})^2\dx
+\int_\R \vr V^2(\px\bar v+1)^2 \dx \\
&\leq 
\left( \|(\px\bar v+1)\|^2_{L^\infty((0,T)\times\R)} + 1\right) 
\int_\R {\vr V^2} \dx 
+\frac{1}{4}\int_\R {\lr{ \frac{\ep\gamma}{\gamma-1}}^2\vr(\px\vr^{\gamma-1})^2} \dx.
\end{split}
\end{equation}
The first term can be controlled by the Gronwall inequality, while the second one can be absorbed by the l.h.s. of \eqref{rel_6}.

\medskip
\noindent{\emph{Estimate of $I_2$}.}  Using the same decomposition of $U$ as in case of $I_1$ we obtain
\eq{\label{I2}
I_2&=-\intR{\vr U\frac{\ep\gamma}{\gamma-1}\partial_x \bar \varrho^{\gamma-1}}\\
&=-\intR{\vr V\frac{\ep\gamma}{\gamma-1}\partial_x \bar \varrho^{\gamma-1}}
+\intO{\vr\lr{ \frac{\ep\gamma}{\gamma-1}}^2
\left(\partial_x \bar \varrho^{\gamma-1}\right) 
\left(\partial_x  \varrho^{\gamma-1}\right)}\\
&\leq \lr{ \frac{\ep\gamma}{\gamma-1}}^2\intR{\lr{\partial_x \bar \varrho^{\gamma-1}}^2\vr}+\frac{1}{4}\intR{\vr V^2}\\
&\quad+
\frac{1}{4}\intR{\lr{ \frac{\ep\gamma}{\gamma-1}}^2\vr(\px\vr^{\gamma-1})^2}+\lr{ \frac{\ep\gamma}{\gamma-1}}^2\intR{(\px\bar\vr^{\gamma-1})^2\vr}\\
&=\frac{1}{4}\intR{\vr V^2}+
\frac{1}{4}\intR{\lr{ \frac{\ep\gamma}{\gamma-1}}^2\vr(\px\vr^{\gamma-1})^2}+2\lr{ \frac{\ep\gamma}{\gamma-1}}^2\intR{(\px\bar\vr^{\gamma-1})^2\vr}.
}
Again, the first two terms will be handled by the Gronwall inequality and by moving to the l.h.s. of \eqref{rel_6}, respectively. The last term on the r.h.s. of \eqref{I2} is a small constant due to the presence of $\ep$, provided that 
$\px\bar\vr^{\gamma-1}\in L^2(0,T; L^2(\R))$. Indeed, from \eqref{rhonorm}, we may deduce that $\ep^2\|\vr\|_{L^\infty(0,T; L^\infty(\R))}\leq C \ep ^{\frac{4(\gamma-1)}{2\gamma-1}}$.

\medskip
\noindent{\emph{Estimate of $I_3$}.}
To deal with this term we again separate $W(x)$ into two parts $K(x)=-|x|$, and $L(x)=\frac{x^2}{2}$
\eq{\label{I3}
I_3=\int_\R{(\vr-\bar\vr)\bar v\px K\ast(\vr-\bar\vr)} \dx +\int_\R{(\vr-\bar\vr)\bar v\px L\ast(\vr-\bar\vr)}\dx.
}
 The part of the potential coming from the 
Newtonian repulsion is good because we can write
\eq{
\partial_{xx}K\ast(\vr-\bar\vr)=-2(\vr-\bar\vr)
}
and so
\begin{equation}\label{Kest_ent}
\begin{split}
\int_\R &
(\vr-\bar\vr)\bar v\px K\ast(\vr-\bar\vr) \dx 
 = -\frac{1}{2}\int_\R\bar v\partial_{xx} K\ast(\vr-\bar\vr)\px K\ast(\vr-\bar\vr) 
\dx
\\
& = - \frac{1}{4}\int_\R 
\bar v\px |\px K\ast(\vr-\bar\vr)|^2 \dx \\
& =\frac{1}{4}\int_\R 
\px\bar v|\px K\ast(\vr-\bar\vr)|^2 \dx 
\leq \frac{1}{4}\| \px \bar v\|_{\infty} 
\int_\R |\px K\ast(\vr-\bar\vr)|^2 \dx.
\end{split}
\end{equation}
Note that $\lim_{|x|\to\infty}\px K\ast(\vr-\bar\vr)=0$ but $\bar v$ may diverge at infinity, so the lack of boundary terms in above requires an extra justification. Note however that
\eq{
\px K\ast(\vr-\bar\vr)=-\int_\R sgn (x-y) (\vr(y)-\bar\vr(y)) \dy
=-2\int_{-\infty}^x(\vr(y)-\bar\vr(y))\dy,
}
and therefore
\eq{
&\lim_{x\to+\infty}\bar v|\px K\ast(\vr-\bar\vr)|^2\leq \lim_{x\to+\infty}|\sqrt{|\bar v|}\px K\ast(\vr-\bar\vr)|^2\\
&\leq\lim_{x\to+\infty}\left|2\sqrt{x}\int_{-\infty}^x(\vr(y)-\bar\vr(y))\dy\right|^2
=4\left(\lim_{x\to+\infty}\frac{\int_{-\infty}^x(\vr(y)-\bar\vr(y))\dy}{x^{-1/2}}\right)^2\\
&=16\left(\lim_{x\to+\infty}\frac{\vr(x)-\bar\vr(x)}{x^{-3/2}}\right)^2=0,
}
where we used d'H\^opital rule and \eqref{second_mom}, subsequently, in the two last equalities.
Coming back to \eqref{Kest_ent} the term on the r.h.s. can be now controlled by the l.h.s. of \eqref{rel_6}, indeed we have
\begin{equation}
\begin{split}
\int_\R (\vr-\bar\vr)K\ast(\vr-\bar\vr) \dx 
& = -\frac{1}{2}
\int_\R \partial_{xx} K\ast(\vr-\bar\vr)K\ast(\vr-\bar\vr)\dx \\
& =\frac{1}{2} \int_\R|\px K\ast(\vr-\bar\vr)|^2 \dx,
\end{split}
\end{equation}
where again we dropped the boundary term due to the fact that $\lim_{|x|\to\infty}\px K\ast(\vr-\bar\vr)=0$.
Therefore, the Newtonian repulsion part of potential $W$  can be controlled by the Gronwall inequality. 
For the quadratic confinement $L(x)=\frac{x^2}{2}$ in the integral in \eqref{I3},  we have by \eqref{Lzero} that
\eq{
\intO{(\vr-\bar\vr)\bar v\left( \px L\ast(\vr-\bar\vr) \right)}
= 0.
}
\medskip
\noindent{\emph{Estimate of $I_4$}.}
The fourth term on the r.h.s. of \eqref{rel_6} can be immediately absorbed by the last term on the l.h.s., we have
\eq{\label{I4}
I_4=2\ep\intR{\vr^\gamma\bar\vr }\leq \ep^{\frac{1}{\gamma+1}}\ep\intO{\vr^{\gamma+1}}+ \ep^{\frac{1}{\gamma+1}}\intO{\bar\vr^{\gamma+1}},}
and so, for sufficiently small $\ep$ the first term is absorbed by $\ep\intO{\vr^{\gamma+1}}$ on the l.h.s. of \eqref{rel_6}, while the second term is arbitrary small provided that $\bar\vr\in L^{\gamma+1}((0,T)\times\R)$.

\medskip
\noindent{\emph{Estimate of $I_5$}.}
Estimate of the last term is almost straightforward as we have
\eq{  I_5 =  -\gamma\intO{\px \bar u F}  \leq  \gamma\|\px\bar v\|_{L^\infty((0,T)\times\R)} \intO{F}}
that can be treated by the Gronwall inequality.

\bigskip

Let us summarise, collecting the estimates for $I_1 - I_5$, we obtain
\eq{\label{rel_8}
&\Dt{\cal E}(\vr,v|\bar\vr,\bar u)(t)+\intR{\vr\lr{\frac{\ep\gamma}{\gamma-1}}^2\lr{\partial_x \varrho^{\gamma-1}}^2}+ 2\ep \intO{\vr^{\gamma+1}} \\
&\leq C( \|(\px\bar v+1)\|^2_{L^\infty((0,T)\times\R)})\intR{\vr V^2}+ \frac{1}{4}\| \px \bar v\|_{L^\infty(0,T\times\R)} \intR{|\px K\ast(\vr-\bar\vr)|^2}\\
& \quad+\|\px\bar v\|^2_{L^\infty((0,T)\times\R)}\gamma \intO{F}+\ep C\\
&\leq C {\cal E}(\vr,v|\bar\vr,\bar v)(t)+\ep C,
}
therefore applying the Gronwall inequality we obtain uniformly in $n$ that
\eq{\label{rel_final}
{\cal E}(\vr_n,v_n|\bar\vr,\bar u)(\tau)+\lr{\frac{\ep\gamma}{\gamma-1}}^2\int_0^\tau\!\!\intR{\vr_n\lr{\partial_x \varrho_n^{\gamma-1}}^2}\,\dt+2\ep\int_0^\tau\!\!\intO{\vr^\gamma_n\bar\vr}\,\dt\\
\leq C_1\lr{{\cal E}(\vr_n,v_n|\bar\vr,\bar u)(0)+\ep C_2}.
}
\section{Conclusion of the proof of Theorem  \ref{Thm:main}}
With the uniform estimates from Section \ref{Sec:3} it is relatively easy to repeat the arguments from \cite{MV07} in order to prove strong  convergence of a subsequence of $(\vr_n,\sqrt{\vr_n}u_n)$ to $(\vr,\sqrt{\vr} u)$ as specified in \eqref{c_delta}. The only modification which perhaps is worth  explaining at this stage is the strong convergence of the density that we now have in much stronger sense, but also only for a restricted range of $\gamma$'s.
\begin{lemma}
Let $\gamma\in(1,\frac32]$, then the sequence $\vr_n^{\gamma-\frac{1}{2}}$ satisfies
\eq{\label{rhon}
\vr_n^{\gamma-\frac{1}{2}}\ is\ bounded \ in \ L^\infty(0,T; H^1(\R));\\
\pt\vr_n^{\gamma-\frac{1}{2}}\ is\ bounded \ in \ L^\infty(0,T; H^{-1}(\R)).
}
As a consequence, up to a subsequence, $\vr_n^{\gamma-\frac{1}{2}}$ converges almost everywhere and strongly in $C([0,T]\times B)$ for any compact subset $B\subset\R$. Moreover $\vr_n$ converges to $\vr$ in 
$C([0,T)\times B)$.
\end{lemma}
\pf The first estimate of \eqref{rhon} follows from \eqref{rhonormg} together with the conservation of mass  \eqref{mass_cons}. Next one can write
$$\pt\vr_n^{\gamma-\frac{1}{2}}=-\px(\vr_n^{\gamma-\frac{1}{2}} u_n)-\lr{\gamma-\frac32}\vr_n^{\gamma-\frac{1}{2}}\px u_n,$$
which yields the second estimate due to \eqref{rhonorm}. The Aubin's Lemma therefore provides the strong local convergence in the space of continuous functions.
To prove the strong convergence of $\vr_n$ it is enough to apply the mean value theorem to the function $f(\vr)=\vr^{\gamma-\frac{1}{2}}$ together with the uniform (in $n$) bound for the density \eqref{rhonorm}. This is the moment of the proof when the assumption $\gamma\leq3/2$ becomes important. $\Box$

The rest of convergences in \eqref{c_delta}, i.e. 
\eq{\label{c_delta1}
&\px\lr{\vr_n^{\gamma-{\frac12}}}\to\px\lr{\vre^{\gamma-{\frac12}}}\quad\text{weakly in } L^2(0,T; L^2(\Omega)),\\
&\sqrt{\vr_n} u_n\to\sqrt{\vre} u_\ep, \quad\text{strongly in } L^{2}(0,T; L^{2}(\R)),\\
&\vr_n v_n\to\vre v_\ep, \quad\text{strongly in } L^{2}(0,T; L^{2}(\R)),\\
&\gamma\vr_n^\gamma\px u_n\to\Lambda_\ep \quad\text{weakly in } L^2(0,T; L^2_{loc}(\R)),
}
follows  the same way as in \cite{MV07} and \cite{HaZa}.

Having proven \eqref{c_delta1}, we can pass to the limit $n\to\infty$ in all terms of system \eqref{main1} to obtain a weak solution $(\vre,\sqrt{\vre}u_\ep)$ as specified in Definition \ref{Def:1}. This is a very similar argument to the one performed in \cite{HaZa}, modulo the passage to the limit in the nonlocal term. Here however, we can use uniform boundedness of higher moments of $\vr_n$, see \eqref{higher_m}, together with strong convergence of $\vr_n$.
This concludes the proof of the first part of Theorem \ref{Thm:main}.

To prove  the second part  we realize  that we can pass to the limit $n\to\infty$ in  \eqref{rel_final}. In this process we basically use the lower semicontinuity of convex functions that together with weak convergence of  $\px \vr_n^{\gamma-\frac12}$ allows us to pass to the limit on the l.h.s. of \eqref{rel_final} to obtain
\eq{\label{rel_final_ep}
{\cal E}(\vr_\ep,v_\ep|\bar\vr,\bar u)(\tau)
\leq C_1\lr{{\cal E}(\vr_\ep,v_\ep|\bar\vr,\bar u)(0)+\ep C_2}.
}
Therefore if $\ep\to 0$, and ${\cal E}(\vre,v_\ep|\bar\vr,\bar u)(0)=0$, then ${\cal E}(\vre,v_\ep|\bar\vr,\bar u)(\tau)\to 0$ for a.e. $\tau\in(0,T)$. $\Box$

\section*{Acknowledgements}
JAC was partially supported by the EPSRC grant EP/P031587/1. AWK is partially supported by a Newton Fellowship of the Royal Society and by the grant Iuventus Plus no. 0871/IP3/2016/74 of Ministry of Sciences and Higher Education RP. EZ was supported by the UCL Department of Mathematics Grant and grant Iuventus Plus  no. 0888/IP3/2016/74 of Ministry of Sciences and Higher Education RP.


\bibliographystyle{plain}

\end{document}